\documentclass{amsart}

\newtheorem{theorem}{Theorem}[section]
\newtheorem{lemma}[theorem]{Lemma}
\newtheorem{proposition}[theorem]{Proposition}

\theoremstyle{definition}
\newtheorem{definition}[theorem]{Definition}

\usepackage{amscd,amssymb}

\def\boxtext#1{%
\vbox{%
\hrule
\hbox{\strut \vrule{} #1 \vrule}%
\hrule
}%
}

\begin{document}

\title[Defining Relations of Noncommutative Trace Algebra]
{Defining Relations of Noncommutative Trace Algebra of Two
$3 \times 3$ Matrices}
\author[Francesca Benanti and Vesselin Drensky]
{Francesca Benanti and Vesselin Drensky}
\address{Dipartimento di Matematica ed Applicazioni, Universit\`a di Palermo,
Via Archirafi 34, 90123 Palermo, Italy}
\email{fbenanti@math.unipa.it}
\address{Institute of Mathematics and Informatics,
Bulgarian Academy of Sciences,
          1113 Sofia, Bulgaria}
\email{drensky@math.bas.bg}
\thanks{The research of the first author was partially supported by MIUR, Italy.}
\thanks{The research of the second author was partially supported by Grant
MM-1106/2001 of the Bulgarian Foundation for Scientific Research.}
\maketitle
\centerline{To the anniversary of Amitai Regev -- mathematician
and friend}
\begin{abstract}
The noncommutative (or mixed) trace algebra $T_{nd}$ is generated by
$d$ generic $n\times n$ matrices and by the algebra $C_{nd}$ generated
by all traces of products of generic matrices, $n,d\geq 2$. It is
known that over a field of characteristic 0 this algebra is a
finitely generated free module over a polynomial subalgebra $S$ of
the center $C_{nd}$. For $n=3$ and $d=2$ we have found explicitly such
a subalgebra $S$ and a set of free generators of the $S$-module
$T_{32}$. We give also a set of defining relations of $T_{32}$ as an algebra
and a Gr\"obner basis of the corresponding ideal. The proofs are
based on easy computer calculations with standard functions of
Maple, the explicit presentation of $C_{32}$ in terms of generators and
relations, and methods of representation theory of the general linear
group.
\end{abstract}


\section*{Introduction}

Let $K$ be any field of characteristic 0 and let
$X_i=\left(x_{pq}^{(i)}\right)$, $p,q=1,\ldots,n$, $i=1,\ldots,d$,
be $d$ generic $n\times n$ matrices. We consider the following two
algebras: the pure (or commutative) trace algebra $C_{nd}$
generated by all traces of products
$\text{\rm tr}(X_{i_1}\cdots X_{i_k})$,
and the mixed (or noncommutative) trace algebra
$T_{nd}$ generated by $X_1,\ldots,X_d$ and $C_{nd}$, where we
treat the elements of $C_{nd}$ as scalar matrices. The algebra
$C_{nd}$ coincides with the algebra of invariants of the general
linear group $GL_n=GL_n(K)$ acting by simultaneous conjugation on
$d$ matrices of size $n\times n$. The algebra $T_{nd}$ is known as
the algebra of matrix concominants and also is the algebra of
invariant functions under a suitable action of $GL_n$. General
results of invariant theory of classical groups imply that the
algebra $C_{nd}$ is finitely generated and $T_{nd}$ is a finitely
generated $C_{nd}$-module. More precise results, see Van den Bergh
\cite{V}, give that $C_{nd}$ is a finitely generated free module
of a polynomial subalgebra $S$. A similar result holds for
$T_{nd}$. Theory of PI-algebras provides upper bounds for the
generating sets of the algebra $C_{nd}$ and of the $C_{nd}$-module
$T_{nd}$. The Nagata-Higman theorem states that the polynomial
identity $x^n=0$ implies the identity $x_1\cdots x_N=0$ for some
$N=N(n)$. If $N$ is minimal with this property, then $C_{nd}$ is
generated by traces of products $\text{\rm tr}(X_{i_1}\cdots
X_{i_k})$ of degree $k\leq N$ and $T_{nd}$ is generated as a
$C_{nd}$-module by products $X_{j_1}\cdots X_{j_l}$ of degree
$l\leq N-1$. These estimates are sharp if $d$ is sufficiently
large. A description of the defining relations of $C_{nd}$ is
given by the Razmyslov-Procesi theory \cite{R, P} in the language
of ideals of the group algebras of symmetric groups. For a
background on the algebras of matrix invariants and concominants
see, e.g. \cite{F2, DF}.

Explicit minimal sets of generators of $C_{nd}$ and $T_{nd}$ and
the defining relations between them are found in few cases only.
It is well known that, in the Nagata-Higman theorem,
$N(2)=3$ and $N(3)=6$, which gives bounds for the degrees of the generators
of the algebars $C_{2d}$ and $C_{3d}$ and their modules $T_{2d}$ and
$T_{3d}$, respectively.

By a theorem of Sibirskii \cite{S}, $C_{2d}$ is generated
by $\text{\rm tr}(X_i)$, $1\leq i\leq d$, $\text{\rm tr}(X_iX_j)$,
$1\leq i\leq j\leq d$,
$\text{\rm tr}(X_iX_jX_k)$, $1\leq i<j<k\leq d$.
There are no relations between
the five generators of the algebra $C_{22}$, i.e.
$C_{22}\cong K[z_1,\ldots,z_5]$. For $d=3$, Sibirskii \cite{S} found one relation
and Formanek \cite{F1} proved that all relations follow from it.
The center of $GL_2$ acts trivially by conjugation, and one has also a
natural action of $PSL_2(K)$. Since $PSL_2({\mathbb C})$ is isomorphic to
$SO_3({\mathbb C})$, one may apply invariant theory of orthogonal groups, see
Procesi \cite{P} and Le Bruyn \cite{L}.
Drensky \cite{D2} translated the generators and defining relations
of the invariants of $SO_3(K)$ and obtained the defining relations
of $C_{2d}$ for all $d$.
As a by-product of a result
of Drensky and Koshlukov \cite{DK} on polynomial identities
of Jordan algebras, one can easily obtain the defining relations of $T_{2d}$,
see the comments in \cite{D2}.

Teranishi \cite{T} found the following system of generators of $C_{32}$:
\[
\begin{array}{c}
\text{\rm tr}(X),\text{\rm tr}(Y),\text{\rm tr}(X^2),\text{\rm tr}(XY),\text{\rm tr}(Y^2), \\
\\
\text{\rm tr}(X^3),\text{\rm tr}(X^2Y),\text{\rm tr}(XY^2),\text{\rm tr}(Y^3),
\text{\rm tr}(X^2Y^2),\text{\rm tr}(X^2Y^2XY),
\end{array}
\]
where $X$, $Y$ are generic $3 \times 3$ matrices. He showed that
the first ten of these generators form a homogeneous system of
parameters of $C_{32}$ and $C_{32}$ is a free module with
generators 1 and $\text{\rm tr}(X^2Y^2XY)$ over the polynomial
algebra on these ten elements. Hence $\text{\rm tr}(X^2Y^2XY)$
satisfies a quadratic equation with coefficients depending on the
other ten generators. The explicit (but very complicated) form of
the equation was found by Nakamoto \cite{N}, over $\mathbb Z$,
with respect to a slightly different system of generators.
Abeasis and Pittaluga \cite{AP} found
a system of generators of $C_{3d}$, for any $d\geq 2$, in terms of
representation theory of the symmetric and general linear groups,
in the spirit of its use in theory of PI-algebras. Teranishi
\cite{T} found also a minimal set of generators and a homogeneous
system of parameters of $C_{42}$. Recently, Aslaksen, Drensky and
Sadikova \cite{ADS} have found another natural set of eleven
generators of the algebra $C_{32}$ and have given the defining
relation with respect to this set. Their relation is much simpler
than that in \cite{N}.

The purpose of the present paper is to find a polynomial
subalgebra $S$ of $C_{32}$ such that both $C_{32}$ and $T_{32}$
are finitely generated free $S$-modules and systems of free
generators of them. We also give a system of generators and
defining relations of $T_{32}$ and find a Gr\"obner basis of
the corresponding ideal with respect to a suitable ordering. The proofs are based
on representation theory of $GL_2$ and use essentially the results
of \cite{ADS} combined with computer calculations with Maple.
Although some of the main results (especially the Gr\"obner basis)
are quite technical, we believe that they may surve as an
``experimental material'' giving some idea for the general picture,
and the developed methods can be successfully
used for further investigation of generic trace algebras.

\section{Preliminaries}

In what follows, we fix $n=3$ and $d=2$ and denote by $X$ and $Y$
the two generic $3\times 3$ matrices.
It is a standard trick to
replace the generic matrices in $T_{nd}$ with generic traceless
matrices. We express $X$ and $Y$ in the form
\[
X=\frac{1}{3}\text{\rm tr}(X)e+x\quad \text{\rm and}\quad
 Y=\frac{1}{3}\text{\rm tr}(Y)e+y,
\]
where $e$ is the identity $3\times 3$ matrix and $x$, $y$ are
generic traceless matrices.
Then
\begin{equation}\label{replacing with traceless matrices}
C_{32}\cong K[\text{\rm tr}(X),\text{\rm tr}(Y)]\otimes_KC_0,
\quad
T_{32}\cong K[\text{\rm tr}(X),\text{\rm tr}(Y)]\otimes_KT_0,
\end{equation}
where the algebra $C_0$ is generated by the traces of products
$\text{\rm tr}(z_1\cdots z_k)$, $z_i=x,y$, $k\leq 6$.
Moreover, $T_0$ is a $C_0$-module generated by the products
$z_1\cdots z_l$, $l\leq 5$.
By well known arguments, as for
``ordinary'' generic matrices, without loss of generality we may
assume that $x$ is a diagonal matrix. Changing the variables
$x_{ii}$ and $y_{ii}$, we may assume that
\begin{equation}\label{matrices}
x=\left( \begin{array}{ccc}
                x_1 & 0 & 0 \\
                0 & x_2 & 0 \\
                0 & 0 & -(x_1+x_2)\\
             \end{array} \right), \
y=\left( \begin{array}{ccc}
                y_{11} & y_{12} & y_{13} \\
                y_{21} & y_{22} & y_{23} \\
                y_{31} & y_{32} & -(y_{11}+y_{22})\\
             \end{array} \right).
\end{equation}
Till the end of the paper we fix the notation $x,y$ for these two
generic traceless matrices. Now we state in detail the result of
\cite{ADS} in the form we need it in our paper. The generators of
$C_0$ in \cite{ADS} are
$\text{\rm tr}(x^2),\text{\rm tr}(xy),\text{\rm tr}(y^2)$,
$\text{\rm tr}(x^3),\text{\rm tr}(x^2y),\text{\rm
tr}(xy^2),\text{\rm tr}(y^3)$, and the other two generators
$\text{\rm tr}(x^2y^2),\text{\rm tr}(x^2y^2xy)$ of the system of
Teranishi \cite{T} are replaced by the elements
\begin{equation}\label{v}
v=\text{\rm tr}(x^2y^2)-\text{\rm tr}(xyxy),
\end{equation}
\begin{equation}\label{w}
w=\text{\rm tr}(x^2y^2xy)-\text{\rm tr}(y^2x^2yx).
\end{equation}
Consider the following
elements of $C_0$:
\[
u=\left| \begin{array}{cc}
          \text{\rm tr}(x^2)& \text{\rm tr}(xy) \\
          \text{\rm tr}(xy) & \text{\rm tr}(y^2) \\
          \end{array} \right|,
\]
\begin{equation}\label{w1 and the others}
 w_1=u^3,  \quad
 w_2=u^2v, \quad
 w_4=uv^2, \quad
 w_7=v^3,
\end{equation}
\begin{equation}\label{w5}
  w_5=v\left| \begin{array}{ccc}
                \text{\rm tr}(x^2)& \text{\rm tr}(xy) & \text{\rm tr}(y^2) \\
                \text{\rm tr}(x^3)& \text{\rm tr}(x^2y) & \text{\rm tr}(xy^2) \\
                \text{\rm tr}(x^2y) & \text{\rm tr}(xy^2) & \text{\rm tr}(y^3) \\
             \end{array} \right| ,
\end{equation}
\begin{equation}\label{w6}
 w_6=\left| \begin{array}{cc}
          \text{\rm tr}(x^3)& \text{\rm tr}(xy^2) \\
          \text{\rm tr}(x^2y) & \text{\rm tr}(y^3) \\
          \end{array} \right|^2
   -4\left| \begin{array}{cc}
          \text{\rm tr}(y^3)& \text{\rm tr}(xy^2) \\
          \text{\rm tr}(xy^2) & \text{\rm tr}(x^2y) \\
          \end{array} \right|\left| \begin{array}{cc}
          \text{\rm tr}(x^3)& \text{\rm tr}(x^2y) \\
          \text{\rm tr}(x^2y) & \text{\rm tr}(xy^2) \\
          \end{array} \right|,
\end{equation}
\begin{equation}\label{w31}
w_3'=u\left| \begin{array}{ccc}
                \text{\rm tr}(x^2)& \text{\rm tr}(xy) & \text{\rm tr}(y^2) \\
                \text{\rm tr}(x^3)& \text{\rm tr}(x^2y) & \text{\rm tr}(xy^2) \\
                \text{\rm tr}(x^2y) & \text{\rm tr}(xy^2) & \text{\rm tr}(y^3) \\
             \end{array} \right| ,
\end{equation}
\begin{equation}\label{w32}
\begin{array}{c}
  w_3''=5[\text{\rm tr}^3(y^2)\text{\rm tr}^2(x^3)
+\text{\rm tr}^3(x^2)\text{\rm tr}^2(y^3)]\\
\\
   -30[\text{\rm tr}^2(y^2)\text{\rm tr}(xy)\text{\rm tr}(x^2y)\text{\rm tr}(x^3)+
        \text{\rm tr}^2(x^2)\text{\rm tr}(xy)\text{\rm tr}(y^3)\text{\rm tr}(xy^2)]\\
\\
  +3\{[4\text{\rm tr}(y^2)\text{\rm tr}^2(xy)+\text{\rm tr}^2(y^2)
            \text{\rm tr}(x^2)][3\text{\rm tr}^2(x^2y)
+2\text{\rm tr}(xy^2)\text{\rm tr}(x^3)]\\
\\
   +[4\text{\rm tr}^2(xy)\text{\rm tr}(x^2)+\text{\rm tr}^2(x^2)\text{\rm tr}(y^2)]
   [3\text{\rm tr}^2(xy^2)+2\text{\rm tr}(x^2y)\text{\rm tr}(y^3)]\}\\
\\
  -2[2\text{\rm tr}^3(xy)+3\text{\rm tr}(x^2)\text{\rm tr}(xy)\text{\rm tr}(y^2)]
  [9\text{\rm tr}(xy^2)\text{\rm tr}(x^2y)+\text{\rm tr}(x^3)\text{\rm tr}(y^3)].
\end{array}
\end{equation}
The element $w_3''$ can be expressed in another
natural way. Recall that the $K$-linear operator $\delta$
of an algebra $R$ is a derivation if
$\delta(f_1f_2)=\delta(f_1)f_2+f_1\delta(f_2)$ for all $f_1,f_2\in R$.
Let $\delta$ be the derivation of $C_0$
which commutes with the trace and is defined by
$\delta(x)=0$, $\delta(y)=x$.
Then
\[
w_3''=\frac{1}{144}\sum_{i=0}^6(-1)^i\delta^i(\text{\rm tr}^3(y^2))
\delta^{6-i}(\text{\rm tr}^2(y^3)).
\]
The main result of Aslaksen, Drensky and Sadikova \cite{ADS}
is that the algebra $C_0$ is generated by
\[
\text{\rm tr}(x^2), \text{\rm tr}(xy), \text{\rm tr}(y^2),
\text{\rm tr}(x^3), \text{\rm tr}(x^2y), \text{\rm tr}(xy^2),
\text{\rm tr}(y^3), v, w%
\]
subject to the defining relation $f=0$, where
\begin{equation}\label{our defining relation}
f=w^2-\left( \frac{1}{27}w_1-\frac{2}{9}w_2+\frac{4}{15}w_3'+\frac{1}{90}w_3''
    +\frac{1}{3}w_4-\frac{2}{3}w_5-\frac{1}{3}w_6-\frac{4}{27}w_7 \right),
\end{equation}
and the elements $v, w, w_1,w_2,w_3',w_3'',w_4,w_5,w_6,w_7$
are given in {\rm (\ref{v}), (\ref{w}), (\ref{w1 and the others}), (\ref{w5}),
(\ref{w6}), (\ref{w31})} and  {\rm (\ref{w32})}.
Let
\begin{equation}\label{the algebra S0}
S_0=K[\text{\rm tr}(x^2),\text{\rm tr}(xy),\text{\rm tr}(y^2),
\text{\rm tr}(x^3), \text{\rm tr}(x^2y), \text{\rm tr}(xy^2),
\text{\rm tr}(y^3), v].
\end{equation}
Then $C_0$ is a free $S_0$-module with basis $\{1,w\}$ and
\[
C_{32}\cong K[\text{\rm tr}(X), \text{\rm tr}(Y)]\otimes_K
S_0[w]/(f).
\]

Recall that the algebras $C_{nd}$ and $T_{nd}$ have a natural multigrading
which takes into account the degrees of the products $X_{j_1}\cdots X_{j_l}$
and of the traces $\text{\rm tr}(X_{i_1}\cdots X_{i_k})$ with respect to each of the
generic matrices $X_1,\ldots,X_d$. The Hilbert series of $C_{nd}$ and $T_{nd}$
are defined as the formal power series
\[
H(C_{nd},t_1,\ldots,t_d)=\sum_{k_i\geq 0}
\dim(C_{nd}^{(k_1,\ldots,k_d)})t_1^{k_1}\cdots t_d^{k_d},
\]
\[
H(T_{nd},t_1,\ldots,t_d)=\sum_{k_i\geq 0}
\dim(T_{nd}^{(k_1,\ldots,k_d)})t_1^{k_1}\cdots t_d^{k_d},
\]
with coefficients equal to the dimensions of the homogeneous
components $C_{nd}^{(k_1,\ldots,k_d)}$ and
$T_{nd}^{(k_1,\ldots,k_d)}$ of degree $(k_1,\ldots,k_d)$,
respectively. The Hilbert series carry a lot of information for
the algebras. In the sequel we shall need the Hilbert series of
$C_{32}$ and $T_{32}$ calculated, respectively, by Teranishi
\cite{T} and Berele and Stembridge \cite{BS}:
\[
H(C_{32},t_1,t_2)=
\frac{1+t_1^3t_2^3}
{(1-t_1)(1-t_2)q_2(t_1,t_2)q_3(t_1,t_2)(1-t_1^2t_2^2)},
\]
\[
H(T_{32},t_1,t_2)=
\frac{1}{(1-t_1)^2(1-t_2)^2
(1-t_1^2)(1-t_2^2)(1-t_1t_2)^2(1-t_1^2t_2)(1-t_1t_2^2)},
\]
where the commuting variables $t_1$ and $t_2$ count, respectively,
the degrees of $X$ and $Y$ and
\begin{equation}\label{the polynomials q2 and q3}
\begin{array}{c}
q_2(t_1,t_2)=(1-t_1^2)(1-t_1t_2)(1-t_2^2),\\
\\
q_3(t_1,t_2)=(1-t_1^3)(1-t_1^2t_2)(1-t_1t_2^2)(1-t_2^3).\\
\end{array}
\end{equation}
Since the Hilbert series of the tensor product is equal to the product
of the Hilbert series of the factors, and
\[
H(K[\text{\rm tr}(X),\text{\rm tr}(Y)],t_1,t_2)=
\frac{1}{(1-t_1)(1-t_2)},
\]
(\ref{replacing with traceless matrices}) implies that
\[
H(C_{32},t_1,t_2)=\frac{H(C_0,t_1,t_2)}{(1-t_1)(1-t_2)},
\quad
H(T_{32},t_1,t_2)=\frac{H(T_0,t_1,t_2)}{(1-t_1)(1-t_2)}.
\]
In this way,
\[
H(C_0,t_1,t_2)=
\frac{1+t_1^3t_2^3}
{q_2(t_1,t_2)q_3(t_1,t_2)(1-t_1^2t_2^2)}.
\]
We rewrite $H(T_0,t_1,t_2)$ in the form
\begin{equation} \label{Berele-series}
H(T_0,t_1,t_2)=
\frac{p(t_1,t_2)}{q_2(t_1,t_2)q_3(t_1,t_2)(1-t_1^2t_2^2)},
\end{equation}
where
\[
p(t_1,t_2)=(1+t_1+t_1^2)(1+t_2+t_2^2)(1+t_1t_2).
\]

Now we summarize the necessary background on representation theory
of $GL_2$. We refer e.g. to \cite{M} for the general facts
and to \cite{D1}  for the applications in the spirit of the problems
considered here.
The irreducible polynomial representations of $GL_2$ are
indexed by partitions $\lambda=(\lambda_1,\lambda_2)$,
$\lambda_1\geq \lambda_2\geq 0$. We denote by $W(\lambda)$ the
corresponding irreducible $GL_2$-module. The group $GL_2$ acts
in the natural way on the two-dimensional vector space
$K\cdot x+K\cdot y$ and this action is extended diagonally on
the free associative algebra $K\langle x,y\rangle$.
As a $GL_2$-module $K\langle x,y\rangle$ is completely reducible and
\[
K\langle x,y\rangle=\sum_{\lambda}d(\lambda)W(\lambda),
\]
where the multiplicity $d(\lambda)$ is equal to the degree of
the corresponding $S_k$-module, $k=\lambda_1+\lambda_2$,
and can be calculated using e.g. the hook formula. In particular,
\[
d(\lambda_1,0)=1,\lambda_1\geq 0,\quad
d(\lambda_1,1)=\lambda_1-1,\lambda_1\geq 1,
\]
\[
d(2,2)=2,\quad d(3,2)=5.
\]
The module $W(\lambda)$ is generated by a
unique, up to a multiplicative constant,
homogeneous element $w_{\lambda}$ of degree $\lambda_1$ and $\lambda_2$
with respect to $x$ and $y$, respectively,
called the highest weight vector of $W(\lambda)$.
It is characterized by
the following property,
see Koshlukov \cite{K} for the version which we need.

\begin{lemma}\label{criterion for hwv}
Let $\delta$ be the derivation of $K\langle x,y\rangle$
defined by $\delta(x)=0$, $\delta(y)=x$.
If $w(x,y) \in K\langle x,y\rangle$
is homogeneous of degree $(\lambda_1,\lambda_2)$,
then $w(x,y)$ is a highest weight vector for some
$W(\lambda_1,\lambda_2)$ if and only if $\delta(w(x,y))=0$.
\end{lemma}

If $W_i$, $i=1,\ldots,m$,
are $m$ isomorphic copies of $W(\lambda)$
and $w_i\in W_i$ are highest weight
vectors, then the highest weight vector of any submodule $W(\lambda)$
of the direct sum $W_1\oplus\cdots\oplus W_m$ has the form
$\xi_1w_1+\cdots+\xi_mw_m$ for some $\xi_i\in K$.
Any $m$ linearly independent highest weight vectors can serve
as a set of generators of the $GL_2$-module $W_1\oplus\cdots\oplus W_m$.
When $m=d(\lambda)$ and
$W_1\oplus\cdots\oplus W_m\subset K\langle x,y\rangle$,
then it is convenient to choose $w_1,\ldots,w_m$ in the following way.
For the standard $\lambda$-tableau

\bigskip
\centerline{\vbox{\hbox{$T_{\sigma}$\hskip0.5truecm$=$\hskip0.5truecm}
\hbox{\phantom{X}}} \hskip0.5truecm \vbox{\offinterlineskip
\hbox{\boxtext{$\sigma(1)$}\boxtext{$\cdots
$}\boxtext{$\sigma(2\lambda_2
-1)$}\boxtext{$\sigma(2\lambda_2+1)$}\boxtext{$\cdots
$}\boxtext{$\sigma(k)$}}
\hbox{\boxtext{$\sigma(2)$}\boxtext{$\cdots
$}\boxtext{$\phantom{-}\sigma(2\lambda_2)\phantom{-}$}}}}
\bigskip

\noindent corresponding to $\sigma\in S_k$, $k=\lambda_1+\lambda_2$,
we associate a highest weight vector $w(T_{\sigma})$ in
$K\langle x,y\rangle$. When $\sigma=\varepsilon$
is the identity of $S_k$ we fix
\[
w(T_{\varepsilon})=(xy-yx)^{\lambda_2}x^{\lambda_1-\lambda_2}
\]
\[
=\sum_{\rho_1,\ldots,\rho_s\in S_2}
\text{\rm sign}(\rho_1\cdots\rho_s)z_{\rho_1(1)}z_{\rho_1(2)}
\cdots z_{\rho_s(1)}z_{\rho_s(2)}x^{\lambda_1-\lambda_2},
\]
where $z_1=x$, $z_2=y$, $s=\lambda_2$. For $\sigma$ arbitrary,
we define $w(T_{\sigma})$ in a similar way, but the skew-symmetries
are in positions $(\sigma(1),\sigma(2)),\ldots,(\sigma(2s-1),\sigma(2s))$
instead of the positions $(1,2),\ldots,(2s-1,2s)$ (and the positions
with fixed $x$ are $\sigma(2s+1),\ldots,\sigma(k)$ instead of
$2s+1,\ldots,k$).
Recall also the Littlewood-Richardson rule, which in the case of $GL_2$ states
\begin{equation}\label{LRR}
W(a+b,b)\otimes_KW(c+d,d)\cong \sum_{s=0}^c W_2(a+b+d+s,b+d+c-s),
\quad a\geq c.
\end{equation}
Finally, if $W$ is a $GL_2$-submodule or a factor module of
$K\langle x,y\rangle$, then $W$ inherits the grading of
$K\langle x,y\rangle$ and its Hilbert series plays the role of
the $GL_2$-character of $W$: If
\[
W\cong \sum_{\lambda}m(\lambda)W(\lambda),
\]
then
\[
H(W,t_1,t_2)=\sum_{\lambda}m(\lambda)S_{\lambda}(t_1,t_2),
\]
where $S_{\lambda}=S_{\lambda}(t_1,t_2)$
is the Schur function associated with
$\lambda$, and the multiplicities $m(\lambda)$ are determined
by $H(W,t_1,t_2)$. In the case of two variables
$S_{(\lambda_1,\lambda_2)}$ has the simple form
\[
S_{(\lambda_1,\lambda_2)}
=(t_1t_2)^{\lambda_2}(t_1^{\lambda_1-\lambda_2}+t_1^{\lambda_1-\lambda_2-1}t_2+\cdots+
t_1t_2^{\lambda_1-\lambda_2-1}+t_2^{\lambda_1-\lambda_2}).
\]

The action of $GL_2$ on $K\langle x,y\rangle$
is inherited by the algebras $C_{32}$, $T_{32}$, $C_0$, and $T_0$.
For example, the elements $v$ and $w$ in (\ref{v}) and (\ref{w})
generate one-dimensional $GL_2$-modules, isomorphic, respectively,
to $W(2,2)$ and $W(3,3)$. As another example,
the polynomial $p(t_1,t_2)$ from
(\ref{Berele-series}) has the form
\begin{equation} \label{the nominator of the Hilbert series for T}
\begin{array}{c}
p(t_1,t_2)=1+S_{(1,0)}+(S_{(2,0)}+S_{(1,1)})\\
\\
+2S_{(2,1)}+(S_{(3,1)}+S_{(2,2)})+S_{(3,2)}+S_{(3,3)}.
\end{array}
\end{equation}
Let $h_k=h_k(t_1,t_2)$ be the homogeneous component of degree $k$ of
the Hilbert series $H(C_0,t_1,t_2)$, i.e.
\[
H(C_0,t_1,t_2)=h_0+h_1+h_2+\cdots.
\]
Direct calculations using the formula
\[
\frac{1}{1-t}=1+t+t^2+\cdots
\]
show that
\begin{equation}\label{the polynomials h}
\begin{array}{c}
h_0=1=S_{(0,0)},\quad h_1=0,\quad h_2=t_1^2+t_1t_2+t_2^2=S_{(2,0)},\\
\\
h_3=t_1^3+t_1^2t_2+t_1t_2^2+t_2^3=S_{(3,0)},\\
\\
h_4=t_1^4+t_1^3t_2+3t_1^2t_2^2+t_1t_2^3+t_2^4=S_{(4,0)}+2S_{(2,2)},\\
\\
h_5=S_{(5,0)}+S_{(4,1)}+S_{(3,2)},\quad
h_6=2S_{(6,0)}+3S_{(4,2)}+S_{(3,3)},\\
\\
h_7=S_{(7,0)}+S_{(6,1)}+3S_{(5,2)}+S_{(4,3)}.\\
\end{array}
\end{equation}

\section{Generators of the $C_{32}$-module $T_{32}$}

In this section we modify for our purposes the idea of
Abeasis and Pittaluga \cite{AP} used in their description of
the generators of $C_{nd}$. We consider the $C_0$-module $T_0$.
It has a system of generators of degree $\leq 5$. Without
loss of generality we may assume that this system is in
the $K$-subalgebra $R_0$ of $T_0$ generated by $x,y$. Let $U_k$
be the $C_0$-submodule of $T_0$ generated by all products
$z_1\cdots z_l$ of degree $l\leq k$, $z_j=x,y$.
Clearly, $U_k$ is also a $GL_2$-submodule of $T_0$.
Let $T_0^{(k+1)}$ be the homogeneous component of degree $k+1$
of $T_0$. Then the intersection $U_k\cap T_0^{(k+1)}$
is a $GL_2$-module and has a complement
$G_{k+1}$ in $T_0^{(k+1)}$,
which is the $GL_2$-module of the ``new'' generators of degree $k+1$.
We may assume that $G_{k+1}$ is a $GL_2$-submodule of $R_0$.
Then the $GL_2$-module of the generators of the $C_0$-module $T_0$
is
\begin{equation}\label{the generators G}
G=G_0\oplus G_1\oplus \cdots\oplus G_5.
\end{equation}

We use the notation $[x,y]=xy-yx$ and $x\circ y=xy+yx$.

\begin{lemma}\label{hvw and relations 31}
For $\lambda=(3,1)$, the following polynomials in $T_0$
are highest weight vectors:
\[
u_1=[x,y]x^2,
\quad
u_2=[x^2,y]x,
\quad
u_3=[x^3,y],
\]
\[
u_4=\text{\rm tr}(x^2)(x\circ y)-2\text{\rm tr}(xy)x^2,
\quad
u_5=\text{\rm tr}(x^2)[x,y],
\quad
u_6=\text{\rm tr}(x^3)y-\text{\rm tr}(x^2y)x.
\]
The elements $u_1,\ldots,u_6$ satisfy the relations
\begin{equation}\label{relation for 31}
2(u_1+u_2)+u_4-u_5+2u_6=0,
\quad
2u_3-u_5=0.
\end{equation}
\end{lemma}

\begin{proof} The elements $u_1,u_2,u_3$ are the highest weight vectors
corresponding, respectively, to the standard tableaux

\bigskip
\centerline{
\vbox{\offinterlineskip
\hbox{\boxtext{1}\boxtext{3}\boxtext{4}}
\hbox{\boxtext{2}}} \hskip1truecm
\vbox{\offinterlineskip
\hbox{\boxtext{1}\boxtext{2}\boxtext{4}}
\hbox{\boxtext{3}}} \hskip1truecm
\vbox{\offinterlineskip
\hbox{\boxtext{1}\boxtext{2}\boxtext{3}}
\hbox{\boxtext{4}}}
}
\bigskip

\noindent
The other elements $u_4,u_5,u_6$ are homogeneous of degree $(3,1)$.
Using Lemma \ref{criterion for hwv}, one checks
that they are also highest weight vectors.
In order to verify the relations,
we have evaluated the expressions in (\ref{relation for 31})
using Maple and assuming that
$x$ and $y$ have the form (\ref{matrices}).
\end{proof}

\begin{lemma}\label{hvw and relations 22}
For $\lambda=(2,2)$, the following polynomials in $T_0$
are highest weight vectors:
\[
u_1=[x,y]^2,
\quad
u_2=x[x,y^2]+y[y,x^2],
\]
\[
u_3=\text{\rm tr}(x^2)y^2-\text{\rm tr}(xy)(x\circ y)
+\text{\rm tr}(y^2)x^2,
\]
\[
u_4=(\text{\rm tr}(x^2)\text{\rm tr}(y^2)-\text{\rm tr}^2(xy))e,
\quad
u_5=(\text{\rm tr}(x^2y^2)-\text{\rm tr}(xyxy))e.
\]
The elements $u_1,\ldots,u_5$ satisfy the relation
\begin{equation}\label{relation for 22}
3w_1-3w_2+3w_3-2w_4+2w_5=0.
\end{equation}
\end{lemma}

\begin{proof} The considerations are similar to those in the
proof of Lemma \ref{hvw and relations 22}.
The elements $u_1,u_2$ are the highest weight vectors
corresponding, respectively, to the standard tableaux

\bigskip
\centerline{
\vbox{\offinterlineskip
\hbox{\boxtext{1}\boxtext{3}}
\hbox{\boxtext{2}\boxtext{4}}} \hskip1truecm
\vbox{\offinterlineskip
\hbox{\boxtext{1}\boxtext{2}}
\hbox{\boxtext{3}\boxtext{4}}}
}
\bigskip

\noindent
and $u_3,u_4,u_5$ are homogeneous of degree $(2,2)$,
satisfying the conditions of
Lemma \ref{criterion for hwv}. Again, the relation
(\ref{relation for 22}) is obtained using Maple.
\end{proof}

\begin{lemma}\label{hvw and relations 41}
For $\lambda=(4,1)$, the following polynomials in $T_0$
are highest weight vectors:
\[
u_1=[x,y]x^3,
\quad
u_2=[x^2,y]x^2,
\quad
u_3=[x^3,y]x,
\quad
u_4=[x^4,y],
\]
\[
u_5=\text{\rm tr}(x^2)[x,y]x,
\quad
u_6=\text{\rm tr}(x^2)[x^2,y],
\]
\[
u_7=\text{\rm tr}(x^3)(x \circ y)-2\text{\rm tr}(x^2y)x^2,
\quad
u_8=\text{\rm tr}(x^3)[x,y],
\]
\[
u_9=\text{\rm tr}(x^2)[\text{\rm tr}(x^2)y-\text{\rm tr}(xy)x],
\quad
u_{10}=[\text{\rm tr}(x^2)\text{\rm tr}(x^2y)
-\text{\rm tr}(xy)\text{\rm tr}(x^3)]e.
\]
The elements $u_1,\ldots,u_{10}$ satisfy the relations
\begin{equation}\label{relation for 41}
\begin{array}{c}
6u_1-3u_5-2u_8=0,\\
\\
6u_2+u_5-2u_6+3u_7-u_8+u_9+2u_{10}=0.\\
\\
2u_3-u_5=0,\\
\\
6u_4-3u_6-2u_8=0.\\
\end{array}
\end{equation}
\end{lemma}

\begin{proof}
Again, $u_1,u_2,u_3,u_4$ are the highest weight vectors
corresponding, respectively, to the standard tableaux

\bigskip
\centerline{
\vbox{\offinterlineskip
\hbox{\boxtext{1}\boxtext{3}\boxtext{4}\boxtext{5}}
\hbox{\boxtext{2}}} \hskip0.5truecm
\vbox{\offinterlineskip
\hbox{\boxtext{1}\boxtext{2}\boxtext{4}\boxtext{5}}
\hbox{\boxtext{3}}} \hskip0.5truecm
\vbox{\offinterlineskip
\hbox{\boxtext{1}\boxtext{2}\boxtext{3}\boxtext{5}}
\hbox{\boxtext{4}}} \hskip0.5truecm
\vbox{\offinterlineskip
\hbox{\boxtext{1}\boxtext{2}\boxtext{3}\boxtext{4}}
\hbox{\boxtext{5}}}
}
\bigskip

\noindent
and $u_5,\ldots,u_{10}$ are found with
Lemma \ref{criterion for hwv}. The relations
(\ref{relation for 41}) are obtained using Maple.
\end{proof}

The proofs of the next two lemmas are similar.

\begin{lemma}\label{hvw and relations 32}
For $\lambda=(3,2)$, the following polynomials in $T_0$
are highest weight vectors:
\[
u_1=[x,y]^2x,
\quad
u_2=(x[x,y^2]+y[y,x^2])x,
\quad
u_3=[x,y][x^2,y],
\]
\[
u_4=[x^2,yx]y+[y^2x,x]x,
\quad
u_5=[x^3,y]y+(y^2x^2-xyxy)x,
\]
\[
u_6=\text{\rm tr}(x^2)[x,y]y-\text{\rm tr}(xy)[x,y]x,
\quad
u_7=\text{\rm tr}(x^2)[x,y^2]-\text{\rm tr}(xy)[x^2,y],
\]
\[
u_8=\text{\rm tr}(x^3)y^2-\text{\rm tr}(x^2y)(x\circ y)
+\text{\rm tr}(xy^2)x^2,
\]
\[
u_9=(\text{\rm tr}(x^2)\text{\rm tr}(y^2)-\text{\rm tr}^2(xy))x,
\quad
u_{10}=(\text{\rm tr}(x^2y^2)-\text{\rm tr}(xyxy))x,
\]
\[
u_{11}=[\text{\rm tr}(x^2)\text{\rm tr}(xy^2)
-2\text{\rm tr}(xy)\text{\rm tr}(x^2y)
+\text{\rm tr}(y^2)\text{\rm tr}(x^3)]e.
\]
They satisfy the relations:
\begin{equation}\label{relation for 32}
\begin{array}{c}
6u_1-6u_3+4u_6-2u_7-6u_8-2u_9+4u_{10}+u_{11}=0,\\
\\
6u_2-6u_3+2u_6+2u_7-6u_8-2u_9-u_{11}=0,\\
\\
6u_4-2u_6+4u_7+2u_{10}+u_{11}=0,\\
\\
6u_5-4u_6+2u_7-2u_{10}-u_{11}=0.\\
\end{array}
\end{equation}
\end{lemma}

In the above lemma the polynomials $u_1,\ldots,u_5$ correspond,
respectively, to the standard tableaux

\bigskip
\centerline{
\vbox{\offinterlineskip
\hbox{\boxtext{1}\boxtext{3}\boxtext{5}}
\hbox{\boxtext{2}\boxtext{4}}} \hskip1truecm
\vbox{\offinterlineskip
\hbox{\boxtext{1}\boxtext{2}\boxtext{5}}
\hbox{\boxtext{3}\boxtext{4}}} \hskip1truecm
\vbox{\offinterlineskip
\hbox{\boxtext{1}\boxtext{3}\boxtext{4}}
\hbox{\boxtext{2}\boxtext{5}}} \hskip1truecm
\vbox{\offinterlineskip
\hbox{\boxtext{1}\boxtext{2}\boxtext{4}}
\hbox{\boxtext{3}\boxtext{5}}} \hskip1truecm
\vbox{\offinterlineskip
\hbox{\boxtext{1}\boxtext{2}\boxtext{3}}
\hbox{\boxtext{4}\boxtext{5}}}
.}
\bigskip

\begin{lemma}\label{hvw and relations 43}
For $\lambda=(4,3)$, the following polynomials in $T_0$
are highest weight vectors:
\[
u_0=[\text{\rm tr}(x^2y^2xy)-\text{\rm tr}(y^2x^2yx)]x,
\]
\[
u_1=[\text{\rm tr}(xy)\text{\rm tr}(y^2)\text{\rm tr}(x^3)
-2\text{\rm tr}^2(xy)\text{\rm tr}(x^2y)
-\text{\rm tr}(x^2)\text{\rm tr}(y^2)\text{\rm tr}(x^2y)
\]
\[
+3\text{\rm tr}(x^2)\text{\rm tr}(xy)\text{\rm tr}(xy^2)
-\text{\rm tr}^2(x^2)\text{\rm tr}(y^3)]e,
\]
\[
u_2=(\text{\rm tr}(x^2)\text{\rm tr}(y^2)
-\text{\rm tr}^2(xy))[\text{\rm tr}(x^2)y-\text{\rm tr}(xy)x],
\]
\[
u_3=(\text{\rm tr}(x^2y^2)-\text{\rm tr}(xyxy))
[\text{\rm tr}(x^2)y-\text{\rm tr}(xy)x],
\]
\[
u_4=\text{\rm tr}(x^3)[\text{\rm tr}(y^3)x
-2\text{\rm tr}(xy^2)y]-\text{\rm tr}(x^2y)[\text{\rm tr}(xy^2)x
-2\text{\rm tr}(x^2y)y],
\]
\[
u_5=\text{\rm tr}(x^2)[\text{\rm tr}(y^3)x^2
-\text{\rm tr}(xy^2)(x\circ y)+\text{\rm tr}(x^2y)y^2]
\]
\[
-\text{\rm tr}(xy)[\text{\rm tr}(xy^2)x^2
-\text{\rm tr}(x^2y)(x\circ y)+\text{\rm tr}(x^3)y^2],
\]
\[
u_6=2[\text{\rm tr}(xy)\text{\rm tr}(xy^2)
-\text{\rm tr}(y^2)\text{\rm tr}(x^2y)]x^2
\]
\[
-[\text{\rm tr}(x^2)\text{\rm tr}(xy^2)
-\text{\rm tr}(y^2)\text{\rm tr}(x^3)](x\circ y)
\]
\[
+2[\text{\rm tr}(x^2)\text{\rm tr}(x^2y)
-\text{\rm tr}(xy)\text{\rm tr}(x^3)]y^2,
\]
\[
u_7=[\text{\rm tr}(x^2)\text{\rm tr}(xy^2)
-2\text{\rm tr}(xy)\text{\rm tr}(x^2y)
+\text{\rm tr}(y^2)\text{\rm tr}(x^3)][x,y],
\]
\[
u_8=[\text{\rm tr}(x^2)\text{\rm tr}(y^2)-\text{\rm tr}^2(xy)][x,y]x,
\]
\[
u_9=[\text{\rm tr}(x^2y^2)-\text{\rm tr}(xyxy)][x,y]x,
\]
\[
u_{10}=[\text{\rm tr}(x^2)\text{\rm tr}(y^2)
-\text{\rm tr}^2(xy)][x^2, y],
\]
\[
u_{11}=[\text{\rm tr}(x^2y^2)-\text{\rm tr}(xyxy)][x^2,y],
\]
\[
u_{12}=\text{\rm tr}(x^3)[x,y]y^2
-\text{\rm tr}(x^2y)[x,y](x\circ y)+\text{\rm tr}(xy^2)[x,y]x^2,
\]
\[
u_{13}=\text{\rm tr}(x^2)[x,y]^2y-\text{\rm tr}(xy)[x,y]^2x,
\]
and satisfy the relation
\begin{equation}\label{relation for 43}
\begin{array}{c}
18u_0-4u_1-2u_2+10u_3+18u_4-6u_5\\
\\
+9u_6+15u_7-4u_8-4u_{10}+12u_{11}-36u_{12}+12u_{13}=0.\\
\end{array}
\end{equation}
\end{lemma}

\begin{proposition}\label{T is a free S-module}
{\rm (i)} The $GL_2$-module $G=G_0\oplus\cdots\oplus G_5$
generates $T_0$ as a $C_0$-module,
where $G_0,G_1,G_2,G_3,G_4,G_5$ are generated
by the sets of elements $\{u_{(0,0)}\}$,
$\{u_{(1,0)}\}$, $\{u_{(2,0)},u_{(1,1)}\}$,
$\{u_{(2,1)}',u_{(2,1)}''\}$,
$\{u_{(3,1)}, u_{(2,2)}\}$, and
$\{u_{(3,2)}\}$, respectively,
and
\[
\begin{array}{c}
u_{(0,0)}=1,\quad
u_{(1,0)}=x,\quad
u_{(2,0)}=x^2,\quad u_{(1,1)}=[x,y],\\
\\
u_{(2,1)}'=[x,y]x,\quad
u_{(2,1)}''=[x^2,y],\\
\\
u_{(3,1)}=[x,y]x^2,\quad u_{(2,2)}=[x,y]^2,
\quad u_{(3,2)}=[x,y]^2x.\\
\end{array}
\]

{\rm (ii)} Let $G_6=Kw$, where $w\in C_0$ is defined in (\ref{w}).
Then $G\oplus G_6$ generates $T_0$ as an $S_0$-module, where $S_0$
is the subalgebra of $C_0$ defined in (\ref{the algebra S0}).
\end{proposition}

\begin{proof}
(i) Let $R_0^{(k)}$ be the homogeneous component of degree $k$
of the $K$-subalgebra $R_0$ of $T_0$ generated by $x,y$.
It is sufficient to show that every irreducible $GL_2$-submodule
$W(\lambda)$ of $R_0^{(k)}$, $k=\lambda_1+\lambda_2\leq 5$, belongs to
$C_0G$. For $k\leq 2$ this is obvious because $R_0^{(k)}=G_k$.
Similarly, any $W(2,1)\subset R_0^{(3)}$ is contained in the
$GL_2$-submodule generated by $u_{(2,1)}',u_{(2,1)}''$.
The submodule $W(3,0)$ of $R_0^{(3)}$ is generated by $x^3$.
We use the equation
\begin{equation}\label{CHtheorem}
x^3-\frac{1}{2}\text{\rm tr}(x^2)x-\frac{1}{3}\text{\rm tr}(x^3)e=0,
\end{equation}
which follows from the Cayley-Hamilton theorem and derive that
$W(3,0)\subset C_0G_1+C_0G_0$. This implies that $x^4$ and $x^5$,
and hence $W(4,0)$ and $W(5,0)$,
also belong to $C_0G$.

Let $\lambda=(3,1)$. In the notation of
Lemma \ref{hvw and relations 31}, the element $u_1=[x,y]x^2$
coincides with $u_{(3,1)}$, hence belongs to $C_0G$. The elements
$u_4,u_5,u_6$ also belong to $C_0G$. Therefore the relations
(\ref{relation for 31}) give that $u_2$ and $u_3$ are in $C_0G$.
Since every $GL_2$-submodule $W(3,1)$ of $R_0^{(4)}$ is contained in the
$GL_2$-module generated by the highest weight vectors $u_1,u_2,u_3$, we
obtain that every $W(3,1)\subset R_0^{(4)}$ is in $C_0G$.

Let $\lambda=(2,2)$. The element $u_1=[x,y]^2$ in
Lemma \ref{hvw and relations 22} is the same as $u_{(2,2)}$
and hence belongs to $C_0G$.
The elements $u_3,u_4,u_5$ also belong to $C_0G$.
The relation (\ref{relation for 22}) gives that $u_2\in C_0G$
and this completes the case $\lambda=(2,2)$ because every
$W(2,2)\subset R_0^{(4)}$ is a submodule of the $GL_2$-module
generated by $u_1$ and $u_2$.

The cases $\lambda=(4,1)$ and $\lambda=(3,2)$ are similar.
In the former case, the relations
(\ref{relation for 41}) in Lemma \ref{hvw and relations 41} give
that all highest weight vectors $u_1,u_2,u_3,u_4$ are
linear combinations of $u_5,\ldots,u_{10}$ and are in $C_0G$.
In the latter case, the element $u_1=[x,y]^2x$ from
Lemma \ref{hvw and relations 32} is equal to $u_{(3,2)}$. Hence
the relations (\ref{relation for 32}) give that
$u_2,u_3,u_4,u_5$ are linear combinations of $u_1=u_{(3,2)}$
and $u_6,\ldots,u_{11}$. In this way, every
$W(3,2)\subset R_0^{(5)}$ is in $C_0G$.

(ii) Since $C_0$ is a free $S_0$-module with basis $\{1,w\}$,
in virtue of (i), it is sufficient to show that
\[
S_0wG\subset \sum_{k=0}^5S_0wR_0^{(k)}\subset S_0G+S_0w.
\]
For $R_0^{(0)}=G_0=K$ we obtain immediately
$S_0wR_0^{(0)}=S_0w\subset S_0G+S_0w$.
The element $w$ generates the one-dimensional
$GL_2$-module $Kw$ isomorphic to $W(3,3)$.
The $GL_2$-module $R_0^{(1)}=G_1\cong W(1,0)$ is generated by
$x$ and has a basis $\{x,y\}$.
The Littlewood-Richardson rule gives that
$W(3,3)\otimes_KW(1,0)\cong W(4,3)$. The module $W(4,3)$
is generated by $wx$ and is spanned by  $\{wx,wy\}$. The
explicit form of the elements $u_1,\ldots,u_{13}$
in Lemma \ref{hvw and relations 43} shows that they belong to
$S_0G$, and the relation (\ref{relation for 43}) gives that
$u_0=wx$ is their linear combination, hence also belongs
to $S_0G$. In this way, $R_0^{(1)}wG_1=S_0wx+S_0wy\subset S_0G$.
For any product $a=z_1\cdots z_k$ of degree $k=1,\ldots,5$,
where $z_j=x,y$, by (i), we know that
$a=a'+a''\in S_0G+S_0wG$, where $a'\in S_0G$ and $a''\in S_0wG$.
Since $T_0$ is a graded vector space, and $S_0G$, $S_0wG$ are its
graded subspaces, the inequality
$\text{\rm deg}(a)\leq 5<6=\text{\rm deg}(w)$ gives that $a''=0$, i.e.
$a\in S_0G$. Also, any product $b=z_1\cdots z_6$, $z_j=x,y$,
belongs to $S_0G+S_0wG$. Since $\text{\rm deg}(b)=6$, and the only
elements of degree 6 in $S_0wG$ are $Kw$, we obtain that
\[
b\in S_0G+Kw\subset \sum_{l=0}^5S_0R_0^{(l)}+Kw.
\]
From $wz_1\in\sum_{l=0}^5S_0R_0^{(l)}$ and
$R_0^{(6)}\subset \sum_{l=0}^5S_0R_0^{(l)}+S_0w$ we derive
\[
wz_1z_2\in\left(\sum_{l=0}^5S_0R_0^{(l)}\right)z_2\subset
\sum_{l=1}^6S_0R_0^{(l)}\subset\sum_{l=0}^5S_0R_0^{(l)}+S_0R_0^{(6)}
\subset \sum_{l=0}^5S_0R_0^{(l)}+S_0w.
\]
Continuing in this way, we obtain that
\[
wz_1\cdots z_k\in \sum_{l=0}^5S_0R_0^{(l)}+S_0w,\quad k\leq 5,
\]
which completes the proof.
\end{proof}

Now we state the main result of this section.

\begin{theorem}\label{T as a free S-module}
Let $T_{32}$ be the mixed trace algebra generated by the
generic $3\times 3$ matrices $X,Y$, and
let $x,y$ be the generic traceless matrices
\[
x=X-\frac{1}{3}\text{\rm tr}(X)e,\quad
y=Y-\frac{1}{3}\text{\rm tr}(Y)e.
\]
Let $G=G_0\oplus G_1\oplus \cdots\oplus G_5$,
where $G_k$ is the $GL_2$-module generated by the elements of degree $k$ among
\[
\begin{array}{c}
u_{(0,0)}=1,\quad
u_{(1,0)}=x,\quad
u_{(2,0)}=x^2,\quad u_{(1,1)}=[x,y],\\
\\
u_{(2,1)}'=[x,y]x,\quad
u_{(2,1)}''=[x^2,y],\\
\\
u_{(3,1)}=[x,y]x^2,\quad u_{(2,2)}=[x,y]^2,
\quad u_{(3,2)}=[x,y]^2x.\\
\end{array}
\]

{\rm (i)} As a $GL_2$-module
\[
\begin{array}{c}
G\cong W(0,0)\oplus W(1,1)\oplus W(2,0)\oplus W(1,1)\\
\\
\oplus 2W(2,1)\oplus W(3,1)\oplus W(2,2)\oplus W(3,2).\\
\end{array}
\]
The vector space $G$ generates $T_{32}$ as a $C_{32}$-module.

{\rm (ii)} Let
\[
v=\text{\rm tr}(x^2y^2)-\text{\rm tr}(xyxy),\quad
w=\text{\rm tr}(x^2y^2xy)-\text{\rm tr}(y^2x^2yx),
\]
\begin{equation}\label{the algebra S}
S=K[\text{\rm tr}(X),\text{\rm tr}(Y),
\text{\rm tr}(x^2),\text{\rm tr}(xy),\text{\rm tr}(y^2),
\text{\rm tr}(x^3), \text{\rm tr}(x^2y), \text{\rm tr}(xy^2),
\text{\rm tr}(y^3), v],
\end{equation}
$G_6=Kw$.
Then $S$ is isomorphic to the polynomial algebra in ten variables and
$T_{32}$ is a free $S$-module. Any basis of the vector space $G\oplus G_6$
serves as a set of free generators of $T_{32}$.
\end{theorem}

\begin{proof} Combining
Proposition \ref{T is a free S-module} (i) and the equality
(\ref{replacing with traceless matrices}) we obtain that $G$ generates
the $C_{32}$-module $T_{32}$. Together with the equality
(\ref{the algebra S0}) this gives that $S$ is isomorphic to
the polynomial algebra in ten variables and $T_{32}$ is an $S$-module
generated by $G\oplus G_6$. Hence, as a graded vector space,
$T_{32}$ is a homomorphic image of the tensor product
$(G\oplus G_6)\otimes_KS$
and $G\oplus G_6$ is a homomorphic image of
\[
\begin{array}{c}
W(0,0)\oplus W(1,1)\oplus W(2,0)\oplus W(1,1)\\
\\
\oplus 2W(2,1)\oplus W(3,1)\oplus W(2,2)\oplus W(3,2).\\
\end{array}
\]
Hence the Hilbert series of $T_{32}$ and $G\otimes_KS$ coincide
if and only if they are isomorphic as graded vector spaces.
The Hilbert series of $W(0,0)\oplus W(1,1)\oplus W(2,0)\oplus W(1,1)
\oplus 2W(2,1)\oplus W(3,1)\oplus W(2,2)\oplus W(3,2)$ is equal to
\[
\begin{array}{c}
1+S_{(1,0)}+(S_{(2,0)}+S_{(1,1)})\\
\\
+2S_{(2,1)}+(S_{(3,1)}+S_{(2,2)})+S_{(3,2)}+S_{(3,3)},
\end{array}
\]
and this is the expression of $p(t_1,t_2)$ from
(\ref{the nominator of the Hilbert series for T}).
The Hilbert series of $S$ is
\[
H(S,t_1,t_2)=\prod_{i=1}^{10}\frac{1}{1-t_1^{a_i}t_2^{b_i}},
\]
where $(a_i,b_i)$ are the degrees of the generators of $S$. Hence
the denominator of $H(S,t_1,t_2)$ is
\[
\begin{array}{c}
(1-t_1)(1-t_2)(1-t_1^2)(1-t_1t_2)(1-t_2^2)\\
\\
(1-t_1^3)(1-t_1^2t_2)(1-t_1t_2^2)(1-t_2^3)(1-t_1^2t_2^2).\\
\end{array}
\]
This shows that both Hilbert series coincide and $T_{32}$
is a free $S$-module generated by any basis of the vector space
$G\oplus G_6$.
\end{proof}

In the end of the section we shall explain how to find the elements $u_i$ in
Lemmas \ref{hvw and relations 31},
\ref{hvw and relations 22},
\ref{hvw and relations 41},
\ref{hvw and relations 32}, and
\ref{hvw and relations 43}, and the relations between them.
We shall illustrate this on Lemma \ref{hvw and relations 31}.
Assume that we already know that the generators (\ref{the generators G})
of degree $\leq 3$ of the $C_0$-module $T_0$ are
\[
G_0=K\cong W(0,0),\quad G_1\cong W(1,0),\quad
G_2=W(2,0)\oplus W(1,1),\quad G_3\cong 2W(2,1).
\]
(Compare with (\ref{the nominator of the Hilbert series for T})!)
Clearly,
\[
T_0^{(4)}=G_4\oplus(C_0^{(1)}G_3+C_0^{(2)}G_2+C_0^{(3)}G_1+C_0^{(4)}G_0)
\]
(the sum in the parentheses may be not direct). Hence $T_0^{(4)}$
is a homomorphic image of the $GL_2$-module
\[
K\langle x,y\rangle^{(4)}\oplus (C_0^{(1)}\otimes_KG_3)
\oplus (C_0^{(2)}\otimes_KG_2)\oplus (C_0^{(3)}\otimes_KG_1)
\oplus (C_0^{(4)}\otimes_KG_0).
\]
Using (\ref{the polynomials h}) we derive
\[
C_0^{(1)}=0,\quad C_0^{(2)}\cong W(2,0),
\quad C_0^{(3)}\cong W(3,0),
\quad C_0^{(4)}\cong W(4,0)\oplus W(2,2).
\]
The Littlewood-Richardson rule (\ref{LRR}) gives that
\[
\begin{array}{c}
C_0^{(2)}\otimes_KG_2\cong W(4,0)\oplus 2W(3,1)\oplus W(2,2),\\
\\
C_0^{(3)}\otimes_KG_1\cong W(4,0)\oplus W(3,1).\\
\end{array}
\]
In this way we prove the existence of the highest weight vectors
$u_4,u_5\in C_0^{(2)}G_2$ and $u_6\in C_0^{(3)}G_1$.
The component $C_0^{(4)}G_0$ does not contain a submodule $W(3,1)$.
From the equality
\[
K\langle x,y\rangle^{(4)}\cong W(4,0)\oplus 3W(3,1)\oplus 2W(2,2)
\]
we obtain $h_1,h_2,h_3$.
The next step is to find the relations (\ref{relation for 31}).
We consider the matrix equation
\begin{equation}\label{the equation for 31}
r_{(3,1)}=\xi_1u_1+\xi_2u_2+\xi_3u_3+\xi_4u_4+\xi_5u_5+\xi_6u_6=0
\end{equation}
with unknowns $\xi_1,\ldots,\xi_6\in K$. The entries
$(r_{(3,1)})_{ab}$ of the $3\times 3$ matrix $r_{(3,1)}$
are polynomials in $x_i$, $y_{ij}$. Each monomial of
$(r_{(3,1)})_{ab}$ is a linear combination of $\xi_1,\ldots,\xi_6$
which gives an ``ordinary'' homogeneous linear equation between
the $\xi$s. We have established that all solutions of
(\ref{the equation for 31}) are linear combinations of
(\ref{relation for 31}). In our exposition we have not used that we know that
(\ref{relation for 31}) give all the solutions and have
derived this with other, ``computer free'' considerations.
The situation is similar for the cases $\lambda=(2,2),(4,1),(3,2),(4,3)$,
when all the relations between the highest weight vectors in $T_0$ follow,
respectively, from (\ref{relation for 22}), (\ref{relation for 41}),
(\ref{relation for 32}), (\ref{relation for 43}).

\section{The Gr\"obner basis}

In this section we give the Gr\"obner basis of the ideal generated by
the defining relations of the algebra $T_{32}$ with respect to
a system of generators and an ordering chosen in a suitable way.
Usually Gr\"obner bases are defined for
ideals of polynomial algebras and free associative algebras with
coefficients from a field, see e.g. \cite{AL} and \cite{U}.
In our case we prefer to state the result
in terms of ideals of free algebras with coefficients from polynomial
algebras. We use some formalism in the spirit of the one considered by
Mikhalev and Zolotykh \cite{MZ}.

Let $U=\{u_1,\ldots,u_p\}$ and $V=\{v_1,\ldots,v_q\}$
be linearly ordered finite sets and let $[U]$ and
$\langle V\rangle$ be the free abelian semigroup
and the free noncommutative semigroup generated by $U$ and $V$,
respectively. The elements of $[U]$ are the ``usual monomials''
$u_1^{k_1}\cdots u_p^{k_p}$ and the elements of
$\langle V\rangle$ are the words $v_{j_1}\cdots v_{j_s}$,
with the ``usual'' multiplication of noncommutative elements,
as in the free associative algebra. We consider the direct product
$[U]\langle V\rangle$ with an arbitrary
total ordering which extends the orderings of $U$ and $V$,
satisfies the descending chain condition, and is a monoid ordering.
The latter means that if $af>bg$ for some $a,b\in[U]$ and
$f,g\in\langle V\rangle$, then $acfh>bcgh$ and $achf>bchg$
for all $c\in[U]$ and $h\in\langle V\rangle$.
We say that $af$ divides $bg$, if there exist
$c\in[U]$ and $h_1,h_2\in\langle V\rangle$ such that
$bg=cah_1fh_2$. We call the elements of $[U]\langle V\rangle$
(generalized) monomials. Every nonzero element of the free associative algebra
\[
K[U]\langle V\rangle
=K[u_1,\ldots,u_p]\langle v_1,\ldots,v_q\rangle
\]
with polynomial coefficients from $K[U]$ is a finite sum of the form
\[
z=\sum_{a_i\in[U]}\sum_{f_j\in\langle V\rangle}\alpha_{ij}a_if_j,\quad
\alpha_{ij}\in K,\quad a_1f_1>a_2f_2>\cdots.
\]
We denote by $z^0$ the leading monomial $a_1f_1$ of $z$.

\begin{definition}
Let $I$ be a two-sided ideal of $K[U]\langle V\rangle$
and let $I^0$ be the set of leading monomials of $I$.
A subset $B$ of $I$ is called a Gr\"obner basis of $I$
(with respect to the fixed total ordering on $[U]\langle V\rangle$)
if for any $z\in I$ there exists a $z_i\in B$ such that
$z^0$ is divisible by $z_i^0$. Equivalently, the set $B^0$
generates the semigroup ideal $I^0$ of $[U]\langle V\rangle$.
\end{definition}

The Gr\"obner basis $B$ has the property that the subset
of $[U]\langle V\rangle$ of all generalized monomials, which are
not divisible by an element of $B^0$, forms a $K$-basis of the factor
algebra $K[U]\langle V\rangle/I$.

Till the end of the paper we fix the set $U$ to consist of
the commuting variables
\begin{equation}\label{commuting variables for presentation of T}
u_{10},u_{01},u_{20},u_{11},u_{02},u_{30},u_{21},u_{12},u_{03},u_{22},
\end{equation}
and $V$ to consist of the noncommuting variables $x_1,y_1,w_{33}$.
Then the algebra $T_{32}$ is a homomorphic image of
$K[U]\langle V\rangle$, under the homomorphism $\pi$ defined by
\begin{equation}\label{T is an image of free algebra}
\begin{array}{c}
\pi:u_{10}\to \text{\rm tr}(X),\quad \pi:u_{01}\to \text{\rm tr}(Y),\\
\\
\pi:u_{ij}\to \text{\rm tr}(x^iy^j),\quad i+j=2,3,\quad
\pi:u_{22}\to v,\\
\\
\pi:x_1\to x,\quad \pi:y_1\to y,\quad \pi:w_{33}\to w,\\
\end{array}
\end{equation}
where the elements $v$ and $w$ are defined, respectively, in
(\ref{v}) and (\ref{w}).

We define an arbitrary total ordering on
$[U]$, with the only restrictions that it satisfies the descending chain
condition and is compatible with the
multiplication. We assume that
$\text{\rm deg}(x_1)=\text{\rm deg}(y_1)=1$ and
$\text{\rm deg}(w_{33})=6$ and order the elements of $\langle V\rangle$
in the deg-lex way: If $f,g\in\langle V\rangle$ and
$\text{\rm deg}(f)>\text{\rm deg}(g)$, then $f>g$
(first by degree). If $\text{\rm deg}(f)=\text{\rm deg}(g)$,
then we order $f$ and $g$
lexicographically assuming that $w_{33}>x_1>y_1$.
Finally, we assume that $af>bg$, $a,b\in[U]$,
$f,g\in\langle V\rangle$, if $f>g$ or, if $f=g$, then $a>b$.
Our purpose is to find
the Gr\"obner basis of the ideal $\text{\rm Ker}(\pi)$ with respect to this
ordering. We need some more relations in $T_{32}$ which have been found and verified
using Maple.

\begin{lemma} The following relations hold in the algebra $T_0$:
\begin{equation}\label{relation of degree 4}
\begin{array}{c}
6(xy)^2-6y^2x^2+3\text{\rm tr}(y^2)x^2-6\text{\rm tr}(xy)xy
+3\text{\rm tr}(x^2)y^2\\
\\
-2(\text{\rm tr}(x^2)\text{\rm tr}(y^2)-\text{\rm tr}^2(xy))e
+2(\text{\rm tr}(x^2y^2)-\text{\rm tr}(xyxy))e=0,\\
\end{array}
\end{equation}
\begin{equation}\label{relation of degree 6}
\begin{array}{c}
36y^2xyx^2-6\text{\rm tr}(y^2)xyx^2
+12\text{\rm tr}(xy)(yx)^2-12\text{\rm tr}(xy)y^2x^2
-6\text{\rm tr}(x^2)y^2xy\\
\\
+12\text{\rm tr}(xy^2)xyx-12\text{\rm tr}(xy^2)yxx
+12\text{\rm tr}(x^2y)yxy-12\text{\rm tr}(x^2y)yyx\\
\\
+(-\text{\rm tr}(x^2)\text{\rm tr}(y^2)
+4\text{\rm tr}^2(xy)
+4(\text{\rm tr}(x^2y^2)-\text{\rm tr}(xyxy)))xy\\
\\
+2(-\text{\rm tr}(x^2)\text{\rm tr}(y^2)-2\text{\rm tr}^2(xy)
+4(\text{\rm tr}(x^2y^2)-\text{\rm tr}(xyxy)))yx\\
\\
+2(\text{\rm tr}(x^2)\text{\rm tr}(y^3)
-6\text{\rm tr}(xy)\text{\rm tr}(xy^2)
+3\text{\rm tr}(y^2)\text{\rm tr}(x^2y))x\\
\\
+2(\text{\rm tr}(y^2)\text{\rm tr}(x^3)
-6\text{\rm tr}(xy)\text{\rm tr}(x^2y)
+3\text{\rm tr}(x^2)\text{\rm tr}(xy^2))y\\
\\
+(2\text{\rm tr}(x^2)\text{\rm tr}(xy)\text{\rm tr}(y^2)
-\text{\rm tr}^3(xy)+\text{\rm tr}(x^3)\text{\rm tr}(y^3)
-3\text{\rm tr}(x^2y)\text{\rm tr}(xy^2)\\
\\
-(\text{\rm tr}(x^2y^2)-\text{\rm tr}(xyxy))\text{\rm tr}(xy)
-3(\text{\rm tr}(x^2y^2xy)-\text{\rm tr}(y^2x^2yx)))e=0.\\
\end{array}
\end{equation}
\end{lemma}

The partial linearizations of the Cayley-Hamilton theorem for
$3\times 3$ traceless matrices (\ref{CHtheorem}) give the following relations
in $T_0$
\begin{equation}\label{relations of degree 3}
\begin{array}{c}
x^2y+xyx+yx^2-\text{\rm tr}(xy)x-\frac{1}{2}\text{\rm tr}(x^2)y-\text{\rm tr}(x^2y)e=0,\\
\\
xy^2+yxy+y^2x-\frac{1}{2}\text{\rm tr}(y^2)x-\text{\rm tr}(xy)y-\text{\rm tr}(xy^2)e=0,\\
\\
y^3-\frac{1}{2}\text{\rm tr}(y^2)y-\frac{1}{3}\text{\rm tr}(y^3)e=0.\\
\end{array}
\end{equation}

We shall need the following easy combinatorial lemma.

\begin{lemma}\label{leading monomials we need}
The only words of degree $\geq 3$ in $\langle x_1,y_1\rangle$
which do not contain as a subword any of the words
\begin{equation}\label{leading monomials for the kernel}
x_1^3,\quad x_1^2y_1,\quad x_1y_1^2,\quad y_1^3,\quad x_1y_1x_1y_1,\quad y_1^2x_1y_1x_1^2
\end{equation}
are
\begin{equation}\label{reduced monomials for the kernel}
\begin{array}{c}
x_1y_1x_1,\quad y_1x_1^2,\quad y_1x_1y_1,\quad y_1^2x_1,\\
\\
x_1y_1x_1^2,\quad y_1x_1y_1x_1,\quad y_1^2x_1^2,\quad y_1^2x_1y_1,\\
\\
y_1x_1y_1x_1^2,\quad y_1^2x_1y_1x_1.\\
\end{array}
\end{equation}
\end{lemma}

\begin{proof}
The words from (\ref{leading monomials for the kernel}) are of degree $\leq 6$
and these from (\ref{reduced monomials for the kernel}) are of degree
$\leq 5$. Hence, it is sufficient to show that all words of degree 3, 4, 5 and 6 which
do not contain (\ref{leading monomials for the kernel}) are those from
(\ref{reduced monomials for the kernel}). This can be done by direct verification.
For example, we give the list of all words of degree 4, writing in parentheses the
subwords from (\ref{leading monomials for the kernel}):
\[
\begin{array}{c}
(x_1^3)x_1,\quad (x_1^3)y_1,\quad (x_1^2y_1)x_1,\quad x_1y_1x_1^2,\quad y_1(x_1^3),\\
\\
(x_1^2y_1)y_1,\quad (x_1y_1x_1y_1),\quad (x_1y_1^2)x_1,\quad
y_1(x_1^2y_1),\quad y_1x_1y_1x_1,\quad y_1^2x_1^2,\\
\\
x_1(y_1^3),\quad y_1(x_1y_1^2),\quad y_1^2x_1y_1,\quad (y_1^3)x_1,\quad (y_1^3)y_1.\\
\end{array}
\]
\end{proof}

The homomorphism $\pi$ from (\ref{T is an image of free algebra})
maps bijectively the polynomial algebra in ten variables $K[U]$ to
$S$, where $S$ is from (\ref{the algebra S}) and $U$ is the set of variables
(\ref{commuting variables for presentation of T}).
We define an isomorphism $\rho:S\to K[U]$ as the inverse of the isomorphism
$\pi\vert_{K[U]}:K[U]\to S$, i.e.
\[
\begin{array}{c}
\rho:\text{\rm tr}(X)\to u_{10},\quad \rho:\text{\rm tr}(Y)\to u_{01},\\
\\
\rho:\text{\rm tr}(x^iy^j)\to u_{ij},\quad i+j=2,3,\quad
\rho:v\to u_{22}.\\
\end{array}
\]
Any element of the subalgebra $SR_0$ of $T_{32}$, where $R_0$ is
generated by $x,y$, has the form
\[
f=\sum a_zz_1\cdots z_l,\quad a_z\in S,\quad z_j=x,y.
\]
It will be convenient, for a fixed presentation of $f$,
to denote by $\rho(f)$ the element
\[
\rho(f)=\sum \rho(a_z)\rho(z_1)\cdots \rho(z_l)
\in K[U]\langle x_1,y_1,w_{33}\rangle,
\]
where $\rho(x)=x_1$, $\rho(y)=y_1$, (and $\rho(e)=1$). Pay attention, that
$\rho(f)$ depends on the concrete form of $f$.

Now we give the elements of $K[U]\langle x_1,y_1,w_{33}\rangle$ which
will be included in the Gr\"obner basis of $\text{\rm Ker}(\pi)$.

The first two elements are
\begin{equation}\label{f1 and f2}
f_1=w_{33}x_1-x_1w_{33},\quad f_2=w_{33}y_1-y_1w_{33}.
\end{equation}
Let $w_1,w_2,w_3',w_3'',w_4,w_5,w_6,w_7$
be the elements from (\ref{w1 and the others}), (\ref{w5}),
(\ref{w6}), (\ref{w31}) and  (\ref{w32}). We use the relation
(\ref{our defining relation}) and construct
\begin{equation}\label{f3}
f_3=w_{33}^2-\rho\left( \frac{1}{27}w_1-\frac{2}{9}w_2+\frac{4}{15}w_3'+\frac{1}{90}w_3''
    +\frac{1}{3}w_4-\frac{2}{3}w_5-\frac{1}{3}w_6-\frac{4}{27}w_7 \right).
\end{equation}
The element $u_0$ from Lemma \ref{hvw and relations 43}
is equal to $xw$. We rewrite
(\ref{relation for 43}) as
\[
f_{43}(x,y)=18xw+\sum_{j=1}^{13}\alpha_ju_j,
\]
\[
\begin{array}{c}
\sum_{j=1}^{13}\alpha_ju_j=-4u_1-2u_2+10u_3+18u_4-6u_5\\
\\
+9u_6+15u_7-4u_8-4u_{10}+12u_{11}-36u_{12}+12u_{13}.\\
\end{array}
\]
Since $w(y,x)=-w(x,y)$, we derive that $f_{43}(y,x)$ has the form
\[
f_{43}(y,x)=-18yw+\sum_{j=1}^{13}\alpha_ju'_j,\quad u'_j=u_j(y,x).
\]
Now we use $f_{43}(x,y)$ and $f_{43}(y,x)$ to define
\begin{equation}\label{f4}
f_4=18x_1w_{33}+\sum_{j=1}^{13}\alpha_ju\rho(u_j),
\end{equation}
\begin{equation}\label{f5}
f_5=-18y_1w_{33}+\sum_{j=1}^{13}\alpha_ju\rho(u'_j).
\end{equation}
The next four equations come from the Cayley-Hamilton theorem
(\ref{CHtheorem}) and its linearizations
(\ref{relations of degree 3})
\begin{equation}\label{f6-f9}
\begin{array}{c}
f_6=x_1^3-\frac{1}{2}u_{20}x_1-\frac{1}{3}u_{30},\\
\\
f_7=x_1^2y_1+x_1y_1x_1+y_1x_1^2-u_{11}x_1-\frac{1}{2}u_{20}y_1-u_{21},\\
\\
f_8=x_1y_1^2+y_1x_1y_1+y_1^2x_1-\frac{1}{2}u_{02}x_1-u_{11}y_1-u_{12},\\
\\
f_9=y_1^3-\frac{1}{2}u_{02}y_1-\frac{1}{3}u_{03}.\\
\end{array}
\end{equation}
Finally, we define $f_{10}$ and $f_{11}$ using the relations
(\ref{relation of degree 4}) and (\ref{relation of degree 6}):
\begin{equation}\label{f10}
f_{10}=6(x_1y_1)^2-6y_1^2x_1^2+3u_{02}x_1^2-6u_{11}x_1y_1
+3u_{20}y_1^2
+2(-u_{20}u_{02}+u_{11}^2+u_{22}),
\end{equation}
\begin{equation}\label{f11}
\begin{array}{c}
f_{11}=36y_1^2x_1y_1x_1^2-6u_{02}x_1y_1x_1^2
+12u_{11}((y_1x_1)^2-y_1^2x_1^2)
-6u_{20}y_1^2x_1y_1\\
\\
+12u_{12}(x_1y_1x_1-y_1x_1^2)
+12u_{21}(y_1x_1y_1-y_1^2x_1)\\
\\
+(-u_{20}u_{02}
+4u_{11}^2
+4u_{22})x_1y_1
+2(-u_{20}u_{02}-2u_{11}^2
+4u_{22})y_1x_1\\
\\
+2(u_{20}u_{03}
-6u_{11}u_{12}
+3u_{02}u_{21})x_1
+2(u_{02}u_{30}
-6u_{11}u_{21}
+3u_{20}u_{12})y_1\\
\\
+2(u_{20}u_{11}u_{02}
-u_{11}^3+u_{30}u_{03}
-3u_{21}u_{12}
-u_{22}u_{11}
-3w_{33}).\\
\end{array}
\end{equation}

The following theorem is the main result of the section.
Pay attention that the Gr\"obner basis which we give
is minimal (no leading monomials of its elements are divisible by each other)
but is not reduced (some of the summands are not in normal form).

\begin{theorem}
The Gr\"obner basis of the kernel of the natural homomorphism
$\pi:K[U]\langle x_1,y_1,w_{33}\rangle\to T_{32}$
from (\ref{T is an image of free algebra}) with respect to the above defined
ordering of $[U]\langle x_1,y_1,w_{33}\rangle$ consists of the polynomials
$f_1$ -- $f_{11}$ from (\ref{f1 and f2}), (\ref{f3}), (\ref{f4}), (\ref{f5}), (\ref{f6-f9}),
(\ref{f10}), and (\ref{f11}).
\end{theorem}

\begin{proof}
The leading monomials of $f_1$ -- $f_{11}$ are
\[
\begin{array}{c}
f_1^0=w_{33}^2,\quad
f_2^0=w_{33}x_1,\quad f_3^0=w_{33}y_1,\quad
f_4^0=x_1w_{33},\quad f_5^0=y_1w_{33},\\
\\
f_6^0=x_1^3,\quad f_7^0=x_1^2y_1,\quad
f_8^0=x_1y_1^2,\quad f_9^0=y_1^3,\\
\\
f_{10}^0=(x_1y_1)^2,\quad
f_{11}^0=y_1^2x_1y_1x_1^2.\\
\end{array}
\]
Now we work modulo the ideal $\text{\rm Ker}(\pi)$.
Let $a\in [U]$, $f\in \langle x_1,y_1,w_{33}\rangle$
be arbitrary monomials. Using the elements
$f_1$ -- $f_{11}$ from (\ref{f1 and f2}), (\ref{f3}), (\ref{f4}), (\ref{f5}), (\ref{f6-f9}),
(\ref{f10}), and (\ref{f11}), we replace $af$ with a linear combination of
monomials which are lower in the ordering of $[U]\langle x_1,y_1,w_{33}\rangle$.
We continue this process until $af$ is presented as a linear combination of
monomials which do not contain any of the leading monomials $f_1^0$ -- $f_{11}^0$,
i.e. these monomials are reduced.
Hence, without loss of generality we may assume that $af$ is reduced.
If $f$ contains $w_{33}$, then, by (\ref{f1 and f2}), (\ref{f4}), and (\ref{f5}),
it cannot contain subwords $w_{33}x_1,w_{33}y_1,x_1w_{33},y_1w_{33}$. By
(\ref{f3}) it cannot contain $w_{33}^2$. Hence, if $f$ contains $w_{33}$, then
$f=w_{33}$. The leading monomials of $f_6$ -- $f_{11}$ are those from
(\ref{leading monomials for the kernel}). Now we use Lemma \ref{leading monomials we need}.
If $af$ is reduced, then $f$ may be any monomial from $\langle x_1,y_1\rangle$
of degree $\leq 2$ or one of the monomials from
(\ref{reduced monomials for the kernel}). The generating function of the set $M$ of
all monomials
from $\langle x_1,y_1,w_{33}\rangle$ in normal form with respect to
$\{f_1,\ldots,f_{11}\}$ (taking into account the degrees of the elements) is
\[
\begin{array}{c}
H(M,t_1,t_2)=1+(t_1+t_2)+(t_1^2+2t_1t_2+t_2^2)+2(t_1^2t_2+t_1t_2^2)\\
\\
+(t_1^3+2t_1^2t_2^2+t_1t_2^3)
+(t_1^3t_2^2+t_1^2t_2^3)+t_1^3t_2^3\\
\\
=1+S_{(1,0)}+(S_{(2,0)}+S_{(1,1)})\\
\\
+2S_{(2,1)}+(S_{(3,1)}+S_{(2,2)})+S_{(3,2)}+S_{(3,3)},
\end{array}
\]
and this is the polynomial $p(t_1,t_2)$ from
(\ref{the nominator of the Hilbert series for T}).
Hence, the Hilbert series of the free $K[U]$-module with basis $M$ is
\[
H(K[U]M,t_1,t_2)=H(K[U],t_1,t_2)H(M,t_1,t_2).
\]
Since the Hilbert series of $K[U]$ is
\[
H(K[U],t_1,t_2)=\frac{1}{(1-t_1)(1-t_2)q_2(t_1,t_2)q_3(t_1,t_2)(1-t_1^2t_2^2)},
\]
where $q_2,q_3$ are defined in (\ref{the polynomials q2 and q3}),
by (\ref{Berele-series}) we obtain that
\[
H(T_{32},t_1,t_2)=H(K[U]M,t_1,t_2).
\]
As a graded vector space, $T_{32}$ is a homomorphic image of
$K[U]M$ and the coincidence of the Hilbert series shows that
$T_{32}\cong K[U]M$, i.e. the $K$-algebra $K[U]\langle x_1,y_1,w_{33}\rangle/\text{\rm Ker}(\pi)$
has a basis consisting of all normal monomials. This implies that $\{f_1,\ldots,f_{11}\}$ is
a Gr\"obner basis of $\text{\rm Ker}(\pi)$.
\end{proof}

\section*{Acknowledgements}

This project was carried out when the second author visited the
University of Palermo with the financial support of the INdAM,
Italy. He is very grateful for the hospitality and the creative
atmosphere during his stay in Palermo.


\begin{thebibliography}{99}

\bibitem{AP}
S. Abeasis, M. Pittaluga, {\it
On a minimal set of generators for the invariants of $3 \times 3$
matrices},
Commun. Algebra {\bf 17} (1989), 487-499.

\bibitem{AL}
W.W. Adams, P. Loustaunau, {\it An Introduction to
Gr\"obner Bases}, Graduate Studies in Math. {\bf 3}, AMS,
Providence, R.I., 1994.

\bibitem{ADS}
H. Aslaksen, V. Drensky, L. Sadikova,
{\it Defining relations of invariants of two $3 \times 3$ matrices},
preprint, http://xxx.lanl.gov/abs/math.RA/0405389.

\bibitem{BS}
A. Berele, J.R. Stembridge, {\it Denominators for the Poincar\'e
series of invariants of small matrices}, Israel J. Math. {\bf 114}
(1999), 157-175.

\bibitem{D1}
V. Drensky, {\it Free Algebras and PI-Algebras},
Springer-Verlag, Singapore, 1999.

\bibitem{D2}
V. Drensky,
{\it Defining relations for the algebra of invariants of
$2\times 2$ matrices},
Algebras and Representation Theory {\bf 6} (2003), 193-214.

\bibitem{DF}
V. Drensky, E. Formanek,
{\it Polynomial Identity Rings}, Advanced Courses in Mathematics,
CRM Barcelona, Birk\-h\"auser, Basel-Boston, 2004.

\bibitem{DK}
V. Drensky, P. Koshlukov,
{\it Weak polynomial identities for a vector space with a symmetric
bilinear form}, in ``Math. and Education in Math.'',
Proc. of 16th Spring Conf. Union of Bulgarian Mathematicians, Sofia,
Publ. House Bulg. Acad. Sci., 1987, 213-219.

\bibitem{F1}
E. Formanek, {\it Invariants and the ring of generic
matrices}, J. Algebra {\bf 89} (1984), 178-223.

\bibitem{F2}
E. Formanek,
{\it The Polynomial Identities and Invariants of $n \times n$
Matrices}, CBMS Regional Conf. Series in Math. {\bf 78},
Published for the Confer. Board of the Math. Sci. Washington DC,
AMS, Providence RI, 1991.

\bibitem{K}
P. Koshlukov,
{\it Polynomial identities for a family
of simple Jordan algebras}, Commun. Algebra {\bf 16} (1988), 1325-1371.

\bibitem{L}
L. Le Bruyn,
{\it Trace rings of generic 2 by 2 matrices},
Memoirs of AMS, {\bf 66}, No. 363, Providence, R.I., 1987.

\bibitem{M}
I.G. Macdonald, {\it Symmetric Functions and Hall
Polynomials}, Oxford Univ. Press (Clarendon),
Oxford, 1979. Second Edition, 1995.

\bibitem{MZ}
A.A. Mikhalev, A.A. Zolotykh,
{\it Standard Gr\"obner-Shirshov bases of free algebras over rings.
I: Free associative algebras},
Int. J. Algebra Comput. {\bf 8} (1998), No. 6, 689-726.

\bibitem{N}
K. Nakamoto,
{\it The structure of the invariant ring of two matrices of degree $3$},
J. Pure Appl. Algebra {\bf 166} (2002), No. 1-2, 125-148.

\bibitem{P}
C. Procesi, {\it The invariant theory of $n\times n$
matrices}, Adv. Math. {\bf 19} (1976), 306-381.

\bibitem{R}
Yu.P. Razmyslov, {\it Trace identities of full matrix
algebras over a field of characteristic zero} (Russian),
Izv. Akad. Nauk SSSR, Ser. Mat. {\bf 38} (1974), 723-756.
Translation: Math. USSR, Izv. {\bf 8} (1974), 727-760.

\bibitem{S}
K.S. Sibirskii,
{\it Algebraic invariants for a set of matrices} (Russian),
Sib. Mat. Zh. {\bf 9} (1968), No. 1, 152-164.
Translation: Siber. Math. J. {\bf 9} (1968), 115-124.

\bibitem{T}
Y. Teranishi,
{\it The ring of invariants of matrices},
Nagoya Math. J. {\bf 104} (1986), 149-161.

\bibitem{U}
V.A. Ufnarovski,
{\it Combinatorial and asymptotic methods
in algebra}, in A.I. Kostrikin, I.R. Shafarevich (Eds.),
``Algebra VI'', Encyclopedia of Math. Sciences {\bf 57},
Springer-Verlag, 1995, 1-196.
\bibitem{V}
M. Van den Bergh, {\it Trace rings are Cohen-Macaulay},
J. Amer. Math. Soc. {\bf 2} (1989), 775-799.

\end{thebibliography}
\end{document}